\documentstyle[12pt]{article}
\normalbaselines
\font\ab=msbm10 scaled\magstep1
\font\aba=msbm7 
\font\ac=eufm10 scaled\magstep1
\newcommand{\oo}{\openone}
\newcommand{\inee}{\ab \symbol{'044}}
\newcommand{\ine}{\mbox{\inee}}

\newcommand{\gz}{\mbox{\ab Z}}
\newcommand{\gr}{\mbox{\ab R}}
\newcommand{\gcp}{\mbox{\ab CP}}
\newcommand{\gc}{\mbox{\ab C}}
\newcommand{\gcm}{\mbox{\aba C}}

\newcommand{\Gras}{\mbox{$G_n({\gc}^{m+n})$}}
\newcommand{\got}[1]{{\mbox{\ac{#1}}}}
\newcommand{\mb}[1]{{\mbox{\boldmath{$#1$}}}}
\newtheorem{rem}{Remark}

\newtheorem{lem}{Lemma}
\newtheorem{th}{Theorem}
\newtheorem{cor}{Corollary}

\catcode `\@=11

\def\openone{\leavevmode\hbox{\small1\kern-3.8pt\normalsize1}}%

\def\@revmess#1#2{\typeout{REVTeX #1: #2}}

\newif\ifsecnumbers \global\secnumbersfalse
\def\ds@eqsecnum{\global\secnumberstrue}

\ifx\selectfont\undefined %
\@revmess{message}{NFSS not detected. Assuming OFSS.}
\fi
\let\reset@font\relax

\newcount\@indentflag \global\@indentflag=1 %
\newdimen\@eqtoeqnum \@eqtoeqnum=6pt %
\def\@indentamount{%
\ifcase\@indentflag 0pt\or\@centering\or0pt plus1fil\fi\relax
}
\def\FL{\global\@indentflag=0 }
\def\FR{\global\@indentflag=2 }

\newcounter{eqletter} 

\def\@eqnnum{\hbox{\reset@font\rm(\theequation)}}
\let\make@eqnnum=\@eqnnum %
\def\eqnum#1{\dec@eqnnum \global\def\make@eqnnum{\reset@font\rm(#1)}%
\def\@currentlabel{#1}%
}
\def\inc@eqnnum{\addtocounter{equation}{1}}
\def\dec@eqnnum{\addtocounter{equation}{-1}}
\@definecounter{equation}%
\@addtoreset{equation}{section} %
\def\theequation@prefix{\arabic{section}.} %

\def\mathletters{%
\inc@eqnnum  \setcounter{eqletter}{0}%
\edef\@currentlabel{\theequation}%
\def\theequation{\theequation@prefix\arabic{equation}.\alph{eqletter}}%
\def\inc@eqnnum{\addtocounter{eqletter}{1}}%
\def\dec@eqnnum{\addtocounter{eqletter}{-1}}%
}

\def\[{\relax
\ifmmode\@badmath\else\par\vskip-\lastskip\vskip\abovedisplayskip\fi
\hbox to\hsize\bgroup
\def\label##1{\@revmess{warning}{\string\label{##1} used
in \string\[\space environment after (\theequation)}}%
\hskip\@totalleftmargin\hskip\@indentamount$\displaystyle
}

\def\]{\relax
\ifmmode
$\hskip\@centering\egroup
\else
\@badmath
\fi
\vskip\belowdisplayskip
\global\@indentflag=1 %
\noindent\ignorespaces
}

\newbox\@testboxa
\newbox\@testboxb

\def\equation{\par\vskip-\lastskip\vskip\abovedisplayskip
\inc@eqnnum\let\@currentlabel=\theequation
\setbox\@testboxa=\hbox\bgroup\hskip\@totalleftmargin\hskip\@indentamount
\hbox\bgroup$\displaystyle
}

\def\endequation{$\egroup\hskip\@centering\egroup %
\setbox\@testboxb=\hbox{\make@eqnnum}%
\bgroup
\@tempdima\wd\@testboxa \advance\@tempdima by\wd\@testboxb
\ifcase\@indentflag
\advance\@tempdima by\@eqtoeqnum
\ifdim\@tempdima<\hsize %
\def\@tempa{0}%
\else
\def\@tempa{1}%
\fi
\or
\advance\@tempdima by2\@eqtoeqnum
\ifdim\@tempdima<\hsize %
\def\@tempa{0}%
\else %
\@tempdima\wd\@testboxa \advance\@tempdima by\wd\@testboxb
\advance\@tempdima by\@eqtoeqnum
\ifdim\@tempdima<\hsize %
\def\@tempa{0}%
\setbox\@testboxa\hbox{\hfill\box\@testboxa\kern\@eqtoeqnum}%
\else
\def\@tempa{1}%
\fi
\fi
\or
\advance\@tempdima by2\@eqtoeqnum
\ifdim\@tempdima<\hsize %
\def\@tempa{0}%
\setbox\@testboxb=\hbox{\kern\@eqtoeqnum\make@eqnnum}%
\else
\def\@tempa{1}%
\fi
\fi
\ifnum\@tempa=0 %
\hbox to\hsize{\unhbox\@testboxa\box\@testboxb}%
\else %
\vbox{\hbox to\hsize{\unhbox\@testboxa}%
\vskip6pt %
\hbox to\hsize{\hfil\box\@testboxb}}%
\fi
\egroup
\global\let\make@eqnnum\@eqnnum %
\vskip\belowdisplayskip\noindent\global\@indentflag=1 \global\@ignoretrue
}

\def\eqnarray{\par\vskip-\lastskip\vskip\abovedisplayskip
\inc@eqnnum\let\@currentlabel=\theequation
\global\@eqnswtrue\m@th
\global\@eqcnt\z@
\tabskip\@totalleftmargin\advance\tabskip by\@indentamount\let\\\@eqncr
\halign to\hsize\bgroup\hskip\@centering
$\displaystyle\tabskip\z@{##{}}$&\global\@eqcnt\@ne
\hfil${{}##{}}$\hfil
&\global\@eqcnt\tw@ $\displaystyle\tabskip\z@{##}$\hfil
\tabskip\@centering \if@eqnsw\phantom{\make@eqnnum\kern\@eqtoeqnum}\fi
&\llap{##}\tabskip\z@\cr}

\def\endeqnarray{%
\@@eqncr\egroup
\vskip\belowdisplayskip\noindent
\dec@eqnnum\global\@indentflag=1
\global\let\make@eqnnum\@eqnnum %
\global\@ignoretrue
}

\def\nonumber{\global\@eqnswfalse
\def\label##1{\@revmess{error}{\string\label{##1} used
with \string\nonumber\space before (\theequation)}}%
}

\def\@@eqncr{\let\@tempa\relax
\ifcase\@eqcnt \def\@tempa{& & &}\or \def\@tempa{& &}%
\else \def\@tempa{&}\fi
\@tempa \if@eqnsw\make@eqnnum\inc@eqnnum\fi
\global\@eqnswtrue\global\@eqcnt\z@\global\let\make@eqnnum=\@eqnnum\cr
}

\@namedef{eqnarray*}{\def\@eqncr{\nonumber\@seqncr}%
\def\label##1{\@revmess{warning}{\string\label{##1} used
in eqnarray* environment before (\theequation)}}%
\eqnarray}
\@addtoreset{equation}{section}

\def\theequation{\arabic{section}.\arabic{equation}}

\def\section{\@startsection {section}{1}{\z@}{-3.5ex plus -1ex minus
     -.2ex}{2.3ex plus .2ex}{\normalsize\bf}}
\def\subsection{\@startsection{subsection}{2}{\z@}{-3.25ex plus -1ex minus
 -.2ex}{1.5ex plus .2ex}{\normalsize\bf}}

\def\@bibitem#1{\item\if@filesw \immediate\write\@auxout
{\string\bibcite{#1}{\the\value{\@listctr}}}\fi\ignorespaces}

\def\@cite#1#2{{#1\if@tempswa , #2\fi}} %
\def\@biblabel#1{$^{#1}$} %

\def\@lbibitem[#1]#2{\item[\@biblabel{#1}]\if@filesw
{\def\protect##1{\string ##1\space}\immediate
\write\@auxout{\string\bibcite{#2}{#1}}}\fi\ignorespaces}

\newdimen\WidestRefLabelThusFar
\global\WidestRefLabelThusFar\z@

\def\bibcite#1#2{\global\@namedef{b@#1}{#2}\@SetMaxRefLabel{#1}}

\def\@SetMaxRefLabel#1{%
{%
\setbox0\hbox{$^{\csname b@#1\endcsname}$}%
\ifdim\wd0>\WidestRefLabelThusFar
\global\WidestRefLabelThusFar=\wd0
\fi
}%
}

\def\@citex[#1]#2{%
\if@filesw\immediate\write\@auxout{\string\citation{#2}}\fi
\leavevmode\unskip$^{\scriptstyle\@cite{\@collapse{#2}}{#1}}$}

\def\CITE{%
\@ifnextchar[{\@tempswatrue\@CITEX}{\@tempswafalse\@CITEX[]}%
}

\let\onlinecite\CITE

\def\@CITEX[#1]#2{%
\if@filesw\immediate\write\@auxout{\string\citation{#2}}\fi
\leavevmode\unskip\ \@cite{\@collapse{#2}}{#1}}

\let\@bylinecite\cite

\def\@collapse#1{%
{%
\let\@temp\relax
\@tempcntb\@MM
\def\@citea{}%
\@for \@citeb:=#1\do{%
\@ifundefined{b@\@citeb}%
{\@temp\@citea{\bf ?}%
\@tempcntb\@MM\let\@temp\relax
\@warning{Citation `\@citeb ' on page \thepage\space undefined}%
}%
{\@tempcnta\@tempcntb \advance\@tempcnta\@ne
\edef\MyTemp{\csname b@\@citeb\endcsname}%
\def\@tempa{\@temptokena=\bgroup}%
\afterassignment\@tempa %
\@tempcntb=0\MyTemp\relax}%
\ifnum\@tempcntb=0\relax%
\@tempcntb=\@MM
\@citea\MyTemp
\let\@temp = \relax
\else %
\edef\@tempd{\number\@tempcntb}%
\ifnum\@tempcnta=\@tempcntb %
\ifx\@temp\relax %
\edef\@temp{\@citea\@tempd}%
\else
\edef\@temp{\hbox{--}\@tempd}%
\fi
\else %
\@temp\@citea\@tempd
\let\@temp\relax
\fi
\fi
}%
\def\@citea{, }%
}%
\@temp %
}%
}%

\def\references{%
\vskip24pt
\hrule width\hsize\relax
\vskip 1.6cm
\list{\@biblabel{\arabic{enumiv}}}%
{\labelwidth\WidestRefLabelThusFar  \labelsep2pt %
\leftmargin\labelwidth %
\advance\leftmargin\labelsep %
\ifdim\baselinestretch pt>1 pt %
\parsep  4pt\relax %
\else   %
\parsep  0pt\relax %
\fi
\itemsep\parsep %
\usecounter{enumiv}%
\let\p@enumiv\@empty
\def\theenumiv{\arabic{enumiv}}%
}%
\let\newblock\relax %
\sloppy\clubpenalty4000\widowpenalty4000
\sfcode`\.=1000\relax
\small
}

\def\endreferences{%
\def\@noitemerr{\@warning{Empty `thebibliography' environment}}%
\endlist     \let\@SetMaxRefLabel\@gobble
}

\def\thebibliography#1{\references}

\catcode `\@=12

\begin{document}

\vspace*{2in}
\noindent
\begin{center}
{\bf COHERENT STATES AND GEODESICS:\\ CUT LOCUS AND CONJUGATE LOCUS}
\vspace{1.3cm}\\
\end{center}
\noindent
\begin{center}
\begin{minipage}{13cm}
 Stefan Berceanu \\
 Institute of Atomic
Physics, Institute of Physics and Nuclear Engineering,
Department of Theoretical Physics, P. O. Box MG-6, Bucharest-Magurele,
Romania;
 E-mail address: Berceanu@Roifa.Bitnet; Berceanu@Roifa.IFA.Ro \\
\makebox[3mm]{ }\\
\end{minipage}
\end{center}

\vspace*{0.5cm}

\begin{abstract}
\noindent
 The  intimate relationship between  coherent states
and geodesics is pointed out. For homogenous manifolds on which the exponential
 from the Lie algebra
to the Lie group equals the geodesic exponential, and in particular for
symmetric spaces,
 it is proved
that the cut locus of the point $0$ is equal to the set of
coherent vectors orthogonal to  $\vert 0>$. A simple method to calculate
the conjugate locus in Hermitian symmetric spaces with significance in
the coherent state approach is presented. The results are illustrated
on the complex Grassmann manifold.

\end{abstract}
\vspace{1cm}
\hspace{1cm}PACS numbers: 02.40.Ma; 03.65-w

\newpage

\section{\hspace{-4mm}.\hspace{2mm}  INTRODUCTION}

\hspace*{0.8cm}   The coherent states\cite{klauder}
 are an excellent interplay of classical
and quantum mechanics.\cite{sbaa}
 The local construction of Perelomov's homogeneous
coherent states\cite{per} was globalized, including the K\"ahlerian
non-homogeneous
manifolds.\cite{raw}
 Simultaneously, the geometric quantization program\cite{kost} furnishes, at
least in principle, a tool towards the quantization program of Dirac on
differentiable manifolds. Actually, using both the same mathematical
objects from complex geometry,\cite{chern} fibre bundles,\cite{hus}
 algebraic topology,\cite{hirz}..., the
coherent state approach and the geometric quantization are deeply
related. In fact, the coherent state approach offers  a straightforward
recipe for
geometric quantization.\cite{cval}

 Interesting problems in both these
already classical fields have not been yet attacked, however. One of them is
{\it the relationship between coherent states and geodesics}.

The starting point of this paper is the Remark 3 in Ref. \onlinecite{sbl} which
 expresses in the language of coherent states the
property, here  called {\it condition A)}, that for symmetric spaces the
 geodesics emanating from  the point
$o$ of the symmetric spaces are given by  the exponential exp from the
Lie algebra to the Lie group and all the
 geodesics are obtained in such a way.
 The aim
of this report is to explore  farther out this relationship.

Firstly, let us remember some notions related to geodesics. Let us fix a point
$p$ of a complete Riemannian manifold $V$ and a geodesic $\gamma$ emanating
 from  $p$.
Then the {\it cut point\cite{kn}} of $p$ along $\gamma$ is the first point on
 $\gamma$ such that, for
any point $r$  beyond $q$ on $\gamma$, there is a shorter geodesic from $p$
 to $r$
different from $\gamma$. A point $q$ is a {\it conjugate point} of $p$  along
$\gamma$  if
there is a 1-parameter family of geodesics from $p$ to $q$ neighbouring
 $\gamma$.
Equivalent and precise definitions are given in Section 3.  Here we only
stress that the importance of the cut loci lies in the fact they inherit
topological
 properties of the manifold $V$. $V$ may be obtained from ${\bf CL}_p$ by
attaching a $n$-dimensional cell via the map ${\rm Exp}:~ CL_p\rightarrow
{\bf CL}_p$ and ${\bf CL}_p$ is a strong deformation retract of $V\setminus
\{p\}$, where $CL$ (${\bf CL}$) denotes the {\it tangent cut locus} (resp.
 the {\it cut locus}).

In this paper it is found
out that for some homogeneous manifolds there is an intimate connection between
 the cut locus ${\bf CL}_0$ of a point on the manifold
 $\widetilde {\bf M}$
corresponding to a fixed coherent vector, say $\vert 0>$, and {\it the polar
divisor} $ \Sigma _0$, i.e. the locus of coherent vectors orthogonal to
$\vert 0>$. The equality
\begin{equation}
{\bf CL}_0 = \Sigma_0\label{basical}
\end{equation}
is proved under a technical condition for the manifold, called {\it condition
B)}. It is stressed that {\it condition B)} imply {\it condition A)} and the
well known case of Riemannian symmetric spaces\cite{wong,won,sak}
 is contained as a particular case.
 Despite the fact that the equality (\ref{basical}) is
proved only for manifolds which verify {\it condition B)},
  this
remark is attractive even from pure mathematical point of view, due to
the lack of methods to characterise the cut locus as an object of global
differential geometry.\cite{kob} In this paper we illustrate the results on
 the case of the complex Grassmann
manifold  \Gras . The cut locus on \Gras~ is well known.\cite{wong,sak}

Another contribution of this article is contained in Theorem 1 which
proposes a calculation with significance in the coherent state approach
of the conjugate locus for Hermitian symmetric
 spaces. The connection with the coherent state approach consists in the fact
that the parameters $Z$ and $B$ which appear in the geodesic exponential
 map $Z=Z(B)$ are  two different parameters
of the coherent states for symmetric spaces. We also illustrate the method
on the case of the complex Grassmann manifold \Gras . However, the situation in
this case is
more complicated than in the case of the cut locus. In fact, there are two
 main contributions in this field.
Wong\cite{won} has announced the expression of the
 conjugate locus in the Grassmann manifold,
while the calculation of the tangent conjugate locus of Sakai\cite{sak}
shows that Wong's result is incomplete. All these problems are largely
 discussed
 elsewhere,\cite{viitor} where another proof of the result of Sakai on the
 tangent conjugate locus is given and a calculation of the conjugate locus
in \Gras ~ using Theorem 1 is also
given. The only new observation\cite{viitor} is a geometrical
 characterisation of the part
of the conjugate locus not found by Wong\cite{won,sak} as consisting of those
points of the \Gras \ which have at least two of the stationary
 angles\cite{jor} with a fixed $n$-plane
equal. We have included in Sec. IV only the notions necessary to
illustrate the results of this paper on the example of
 \Gras .

Some of the results included in the present work have been already briefly
 announced
as part of a trial to find a geometrical characterisation of Perelomov's
construction of coherent state manifold as K\"ahlerian embedding into a
projective space.\cite{sb,sbpol}

The paper is organised as follows. In Section 2 the notation on
coherent state manifolds is fixed. The result ${\bf CL}_0 = \Sigma_0$ and
some results on coherent states and conjugate points are proved in
Section 3 for  manifolds $\widetilde {\bf M}$ which verify {\it condition B)}.
 Section 4 deals with the Grassmann manifold.

\section{\hspace{-4mm}.\hspace{2mm}   THE COHERENT STATE MANIFOLD AND
THE COHERENT VECTOR MANIFOLD}
\hspace*{0.8cm} Firstly we fix the notation referring to the
coherent state manifold.

1. Let us consider a quantum system with symmetry, i.e. a triplet
$({\bf K},G,\pi)$, where  $\pi$ is an unitary
irreducible representation of the Lie group $G$ on the Hilbert space
${\bf K}$. Let us consider the orbit

\begin{equation}
\widetilde {\bf M} =\{\widetilde \pi(g)\vert
\widetilde\psi_0>~ \vert ~ g
\in G\},
\end{equation}
where  $\tilde \pi$ is the projective representation of $G$ induced
by $\pi,~ \vert \psi_0> \in {\bf K} $ is fixed and $\xi:
{\bf K}\rightarrow
{\bf PK}$ is the projection $\xi(\vert \psi>)\equiv \vert
\widetilde\psi>=\{e^{i\varphi}  \vert \psi>\vert \varphi \in {\gr}\}$.
Then we have the bijection  $\tilde \xi :G/K\rightarrow
\widetilde {\bf M},~ \tilde\xi(gK)=\widetilde \pi(g){\vert \widetilde
\psi_0>}$, where $K$ is the stationary group of  the state
 $\vert \widetilde \psi_0>$. The quantum mechanics  can be realised
as the elementary $G$-space\cite{gs} $({\bf
PK},\omega_{FS},\rho'$), where $\omega_{FS}$ is the
Fubini-Study  (K\"ahler) fundamental two-form on the projective space
 ${\bf PK}$,
and $\rho'$ is the isomorphism of the Lie algebra ${\got g}$
of $G$ into the algebra of smooth functions on ${\bf PK}$.

The keystone in the coherent state approach is to find a Hilbert space
${\bf L}$ and an K\"ahlerian embedding $\iota :\widetilde {\bf M}
\hookrightarrow {\bf PL}.$\cite{sb,od} Then
$\widetilde {\bf M}$ is called {\it
coherent state manifold} and  $(\widetilde {\bf M},\omega,\rho)$ is a
hamiltonian $G$-space, with $\omega= \omega_{FS\vert \widetilde
{\bf M}}=\iota ^*\omega_{FS},~
\rho = \rho^\prime_{\vert \widetilde {\bf M}}$. {\it Dequantization}
means passing on from the dynamical system problem in the initial
Hilbert space {\bf K} to the corresponding one on $\widetilde {\bf M}$.

If $\vert \widetilde \psi_0>\equiv \vert j>$, i.e. an (anti-)dominant
weight vector for compact connected simply connected Lie groups, then
$\iota$ is indeed a K\"ahlerian embedding\cite{gs} and
 $\widetilde {\bf M}$ coincides
with the coadjoint orbit in ${\got g}^*$ through the root $j$  corresponding
 to the (anti-)dominant weight vector.\cite{kir}
 So, furnishing both the representation
$\pi = \pi_j$ and the Hilbert space $\bf K_j$ of holomorphic sections
with base  $\widetilde {\bf M}$, it is found out that ${\bf L} ={\bf K}^*_j,$
 and the Borel-Weil-Bott theorem solves the requantization
problem.\cite{cval,on} Here
$E^*$ denotes the dual of the vectorial space $E$,
 i.e. the space of linear functionals on $E$.

Now we briefly discuss the embedding $\iota$ for compact complex manifolds
 $\widetilde {\bf M.}$  In this case, {\it the condition for the existence
of the embedding} $\iota$ {\it is equivalent with the requirement for the
manifold to be Hodge,}\cite{hirz}
{\it which is the same condition as prequantization in
geometric quantization.}\cite{wood}
 For example, in order to have the condition
 $\omega \in H^2( \widetilde {\bf M},\gz)$ fullfiled,   for Hermitian
 symmetric spaces
 it is sufficient
 that a  theorem due to Harish-Chandra\cite{sbl,kost,knapp}
to be satisfied. This theorem
  in the compact case is just the Borel-Weil-Bott theorem.
 The Kodaira vanishing theorem replaces the Borel-Weil-Bott
theorem, as was already remarked\cite{cgr} in the context of cohernt states.
Let now $\xi _0:{\bf M}' \rightarrow \widetilde {\bf M}$ be a holomorphic
 line bundle.
Another way to express the condition to have the embedding
 is that the line bundle $\bf M'$ be a
positive one, or, equivalently, to be ample (see Thm. 5.1.
p. 89 in Ref. \onlinecite{ss}).
 The last condition
means that there exists an integer $m_0$ such that for $m\geq m_0,~
M\equiv M'^m=\iota ^*[1]$. We use the notation $[r]=H^r,\mbox{} r\in {\gz}$,
 where $H$ is
the hyperplane bundle over ${\bf PL}$ and $E^m$ is the $m-$times tensor
product of the
bundle $E$ with itself. Here $\xi _0$ is the positive line bundle appearing
in the Kodaira embedding theorem, and the embedding $\iota :\widetilde {\bf M}
\hookrightarrow {\bf PL}=\gcp^{N-1}$  is\cite{gh}

\begin{equation}
 \iota\equiv \iota_{\bf M}: x\rightarrow \iota_{\bf M}(x)=[s_1(x),...,s_N(x)] .
\end{equation}

The line bundle ${\bf M}$ is furnished
by the coherent state approach and is called
{\it coherent vector manifold}.\cite{sbcag}
As a consequence of the Kodaira embedding theorem, the
Kodaira
vanishing\cite{ss} theorem implies that in the sum giving the generalised
 Euler-Poincar\'e characteristic,\cite{hirz}
 only the zero term is present, and the
dimension of the representation $\pi _j$ is furnished by the
 Riemann-Roch-Hirzebruch theorem (cf. Thm. 18.2.2 p. 140 in Ref.
\onlinecite{hirz}). There are situations in which the
coherent state approach permits rapid and explicit statements, for example,
for flag manifolds,
the minimal exponent $N$ appearing in the Kodaira embedding theorem,
$\widetilde {\bf M}\hookrightarrow {\bf CP}^{N-1}$, is equal
with the Euler-Poincar\'e characteristic,
 $N=\chi(\widetilde {\bf M})$.\cite{sb,sbcag}

The noncompact case is treated similarly by Kobayashi,\cite{kobi}
 the Hilbert
space {\bf L} being infinite dimensional. In  the construction of Kobayashi,
$\bf L$ is the dual of the Hilbert space of square integrable holomorphic
$n$-forms in $\widetilde {\bf M}$. If $K$ is the kernel $2n$-form on
$\widetilde {\bf M}\times \overline{\widetilde {\bf M}}$, then the K\"ahler
metric used by Kobayashi is $ds^2=\sum \partial^2\log K^*/\partial z_i\partial
\overline{z}_j$, where $K(z,\overline{z})=K^*(z,\overline{z})dz_1\wedge
\ldots\wedge dz_n\wedge d\overline{z}_1\wedge\ldots\wedge d\overline{z}_n$.

 The condition A1) (A2),
respectively A3)) in Kobayashi corresponds to the condition of the set
of divisors  without base points (the differential of $\iota$ do
not has degenerate points,  respectively, the condition A1) plus
the injectivity condition in the book of Griffith and Harris\cite{gh}). We
remember that A1) implies that $ \omega_{\vert
 \widetilde{\bf M}}= \iota^*\omega_{FS}$, while A2) and A3) implies that
the application $\iota$ is a K\"ahlerian embedding.

Now we discuss other cases in which the representation $\pi$
can be constructed. The condition to have holomorphic discrete series
on homogeneous bounded symmetric domains (non-compact Hermitian
symmetric spaces) results from the quoted theorem of Harish-Chandra
and, more generally, the condition to have  discrete series for
connected semisimple Lie groups is that
 ${\rm rank}\, G= {\rm rank}\, K$.\cite{knapp}
 The problem
of structure of homogeneous K\"ahler manifolds in the context of fundamental
conjecture has began to be handled in connection with the coherent
states, especially for the unimodular groups.\cite{lis}

2. In this section we restrict ourselves to coherent state manifolds of flag
type, i.e. $\widetilde {\bf M} \approx G/K \approx G^{\gcm}/P$, where $G$
is a compact connected simply connected semisimple Lie group, $G^{\gcm}$
 is
the complexification of $G$ and $P$ is a parabolic subgroup of
 $G^{\gcm}$.\cite{sbl}\linebreak The
noncompact case is handled similarly, whenever the conditions of the
existence of the representation  $\pi_j$ are fulfilled.

 Let $W(G)=N(T)/C(T)$ denote the Weyl group associated with $G$, where
$N(T)$ $(C(T))$ is the normalizer (the centralizer) of the
Cartan group T. Let $\Sigma \subset N(T) $ be a set of elements such
that quotient space $W(G)/W(K)$ is made of the coset classes
$\{sC(T)\}W(K), s\in \Sigma$. Then there is an open covering of
  $\widetilde {\bf M}$ by $({\cal V}_s)_{s\in \it \Sigma}$, where ${\cal V}_s=
\pi_j(s){\cal V}_0, s\in \Sigma$.\cite{sbcag} The coherent state vectors
 corresponding
to the points of the neighbourhood ${\cal V}_0 \subset \widetilde {\bf M}$
 around $Z=0$ are

\begin{equation}
\vert
Z,j>=\exp\sum_{{\varphi}\in\Delta^+_n}(Z_{\varphi}F^+_{\varphi})\vert j> ,
\vert {\underline Z}>=<Z\vert Z>^{-1/2}\vert Z> \in  {\bf M} ~,\label{z}
\end{equation}
 where $Z\in {\gc}^n$ are local coordinates and $n$ is the dimension of
the manifold   $\widetilde {\bf M}$.   Here

\begin{equation}
F^\pm_\varphi=\pi^{*\prime}(f^\pm_\varphi),
\end{equation}
$\pi'=d\pi,~ \pi '$ is the isomorphism of the Lie algebra {\got g} of $G$
onto the Lie algebra of operators on {\bf K}, $\pi^*$ is the group
isomorphism
$G^{\gcm} \rightarrow \pi^{*}(G^{\gcm}),$

\begin{equation}
\pi^*(e^Z)=\exp(\pi^{*\prime}(Z)),~ Z\in {\got g}^{\gcm},
\end{equation}
$\pi^{*\prime}({\got g}^{\gcm})$ is the complexification of the Lie algebra
$\pi'({\got g})$, the subindex $n~ (c)$ abbreviates the noncompact
(respectively, compact) , $\Delta$ are the roots and $\Delta ^+$ the
 positive roots.

We also use the notation

\begin{equation}
f^\pm_\varphi=\cases {k^\pm_\varphi=ie_{\pm\varphi},& for $
X_n$,\cr ~~\cr e^\pm_\varphi=e_{\pm\varphi}, &for $X_c$.\cr}
\end{equation}
where $e^\pm_\varphi=e_{\pm\varphi}$ are the part of the Cartan-Weyl base
corresponding to {\got m}. Here ${\got g}={\got k}\oplus{\got m}$ is the Cartan
decomposition of the Lie algebra {\got g} of {\bf G} and {\got k} is the Lie
algebra of {\bf K}.

The homogeneous symmetric spaces  are obtained as

\begin{equation}
X_{n,c}=\exp\sum_{{\varphi}\in\Delta^+_n}(B_{\varphi}f^+_{\varphi}-{\bar
B}_{\varphi}f^-_{\varphi}) \,\cdot~o ,
\end{equation}
where $o=\lambda (e)$, $e$ is the unit element in $G$ and $\lambda$ is
the canonical projection $\lambda :G\rightarrow G/K$.
Let also the notation

\begin{equation}
\vert
B,j>=\exp\sum_{{\varphi}\in\Delta^+_n}(B_{\varphi}F^+_{\varphi}-{\bar
B}_{\varphi}F^-_{\varphi})\vert j> , \label{b}
\end{equation}
\begin{equation}
\vert
B,j>\equiv \vert {\underline {Z,j}}>.
\end{equation}

Note that
\begin{equation}
F^+_{\varphi}\vert j>\not= 0,  ~
F^-_{\varphi}\vert j> = 0,~ H_i\vert j>=j_i\vert j>,\label{act}
\end{equation}
where $ \varphi \in\Delta^+_n$,  $H_i=\pi^*(h_i)$,  $\{h_i\}$ is a
base of the Cartan subalgebra and $ i=1,\ldots,{\rm rank}\,G$.

3. We now state more precisely the definition of the coherent vector
manifold ${\bf M}$ corresponding to flag manifolds $\widetilde {\bf M}$. Let
$ {\bf M'}$ be the holomorphic line bundle
 $ {\bf M'}= \xi ^{-1}_0(\widetilde {\bf M})\rightarrow G_c/P$
 associated by the holomorphic character $\chi =
\chi _j$ of $P$ to the principal bundle $P\rightarrow G^{\gcm}
\rightarrow  G^{\gcm}/P$, i.e. the
line bundle obtained identifying $(g,\chi (p)w)$ with $(gp,w)$, where
$p\in P, w\in {\gc}$.

In fact, if ($ {\bf M'},\omega ,J$) is the compact
 K\"ahler manifold ($J$ is the complex structure, $J=ad(Z)\vert _{\got m}$
and $Z$ is the central element of the Lie algebra {\got k}),
 then ($ {\bf M'},\nabla,h$) is a quantization bundle over
$\widetilde {\bf M}$,\cite{kost} where $h$
is the hermitian form on the tautological line bundle [-1] over ${\bf PL}$.
Then, on [-1], $h$ is given  by $h: z\rightarrow \vert z\vert ^2$. Also
$\mbox\rm{curv}(\nabla )=-2\pi i\omega$, so
$\omega\in c_1({\bf M}')=[\omega ]_{\mbox {de Rham}}$.

 If
$\varphi _i: {\it V}_i\times {\bf C}\rightarrow \xi^{-1}_0({\it V}_i)$ is
the local trivialization of the holomorphic line bundle
${\bf M'}\rightarrow \widetilde {\bf M}$, then a
global section is given by

\begin{equation}
\vert s_i(m)>=(g_i(Z_i),f_{s_i}(Z_i))=(g_i(Z_i),<s_i\vert Z_i>),
\end{equation}
where $m=g_i(Z_i)\in {\it V}_i$ are matrix elements determined by the local
coordinates $Z_i$. Then the scalar product on the line bundle
$ {\bf M'}\rightarrow \widetilde {\bf M}$ is given by\cite{cval,sbcag}

\begin{eqnarray}
\nonumber <s_i\vert s'_i> & = & \int _{\widetilde {\bf M}}h_X(s_i(X),s'_i(X))
{\omega ^n(X)
\over n!}\\
 & = & <f_{s_i},f_{s'_i}>\mbox{}=\int _{\widetilde {\bf M}}h_X(f_{s_i}(X),
f_{s'_i}(X)){d\mu (X)\over<X\vert X>}~,
\end{eqnarray}
where $d\mu (X)$ is the Haar measure on $\widetilde {\bf M}\approx G^{\gcm}/P$.

The scalar product in (2.12) is also a hermitian scalar product of sections
with base $\widetilde {\bf M}$ in the $D_{\widetilde {\bf M}}$- module
of differentiable operators on $\widetilde {\bf M}$.\cite{sbl}

When both the dequantization and the requantization can be done, then
the Hilbert space ${\bf K}_j$ attached to the representation $\pi _j$
and the initial {\bf K} are isomorphic.\cite{cval,cgr}

\section{\hspace{-4mm}.\hspace{2mm}   THE CUT LOCUS AND COHERENT STATES}

\hspace*{0.8cm} In this section we shall be concerned with various
aspects of the relationship between geodesics and coherent states. We
briefly review some definitions used in the Introduction.

1. Let $V$ be compact Riemannian manifold of dimension $n$, $p\in V$ and let
 ${\rm Exp}_p$
be the (geodesic) exponential map at the point $p$. Let $C_p$ denote
the set of vectors $X\in V_p$ (the tangent space at $p\in V$) for which
${\rm Exp}_pX$ is singular. A point $q$ in $V$ ($V_p$) {\it is conjugate
to p} if it is in ${\bf C}_p={\rm Exp}\,C_p$ ($C_p$)\cite{helg} and ${\bf
C}_p$ ($C_p$) is called {\it the conjugate locus} (resp. {\it tangent
 conjugate locus}) of the point $p$.

Let $q\in V$. The point $q$ is in the {\it cut locus} ${\bf CL}_p$ of $p\in V $
if it is nearest point to $p\in V$ on the geodesic joining $p$ with $q$,
beyond which the geodesic ceases to minimise its arc length.\cite{kn}  More
 precisely, let $\gamma _X(t)={\rm Exp}\, tX$ be a geodesic emanating from
 $\gamma_X(0)=p\in V$, where  $X$ is a unit vector from the unit sphere $S_p$
 in $V_p$.
$t_0X$ (resp. ${\rm Exp}\,t_0X$) is called a {\it tangential cut point}
 ({\it cut point}) of $p$ along $t\rightarrow{\rm Exp}\,tX$
 ($0\leq t\leq s$) if
the geodesic segment joining $\gamma_X(0)$ and $\gamma_X(t)$ is a minimal
geodesic for any $s\leq t_0$ but not for any $s>t_0$.

 Let us define the
 function $\mu :
S_p\rightarrow {\gr}^+\cup \infty,~ \mu (X)=r$, if $q={\rm Exp}\,rX\in
{\bf CL}_p$, and $\mu (X)=\infty$ if there is no cut point of $p$ along
$\gamma_X(t)$. Setting $ I_p=\{tX,~0<t<\mu (X)\}$, then ${\bf I}_p=
{\rm Exp}\,I_p$ is called the {\it  interior set at p}. Then:

1) ${\bf I}_p\cap {\bf CL}_p=\emptyset ,~ V={\bf I}_p\cup {\bf CL}_p,$ the
 closure $\bar {\bf I}_p=V$, and dim ${\bf CL}_p\leq n-1$;

2) $I_p$ is a maximal domain  containing $0=0_p\in V_p$ on which
 ${\rm Exp}_p$ is a
 diffeomorphism and ${\bf I}_p$ is the largest open subset of $V$ on which a
normal coordinate system around $p$ can be defined.

The relative position of ${\bf CL}_0$ and ${\bf C}_0$ is given by Theorem 7.1
 p. 97 in Ref. \onlinecite{kn} reproduced below.

 Let the notation
 $\gamma_t=\gamma _X(t)$. Let $\gamma_r$ be the cut point of $\gamma_0$
 along a geodesic $\gamma = \gamma_t,~ 0\leq
t< \infty$. Then, at least one (possibly both) of the following
statements holds:

(1) $\gamma_r$ is the first conjugate point of $\gamma_0$ along $\gamma$;

(2) there exists, at least, two minimising geodesics from $\gamma_0$ to
 $\gamma_r$.

Crittenden\cite{cr} has shown that for the case of simply connected symmetric
spaces, the cut locus is identified to the first conjugate point.  Generally,
 the situation  is more complicated.\cite{war,wei}

Here are simple examples of cut loci. For the sphere $S^n$, the cut locus of a
point reduces to the antipodal point, while the tangent cut locus $CL$ is the
sphere of radius $\pi$
 with centre at the origin of the tangent space. For ${\gcp}^n$, $CL$ is
 also
the sphere of radius $\pi$ with centre at the origin of the tangent space to
${\gcp}^n$ at the given point, while ${\bf CL}$ is the hyperplane at infinity
${\gcp}^{n-1}$. Except few situations, e. g.  the ellipsoid, even for low
 dimensional manifolds as the (asymmetric Berger's spheres) $S^3$, ${\bf CL}$
is not known explicitly. Helgason\cite{helg} has shown that the  cut locus of
a
compact connected Lie group, endowed with a bi-invariant Riemannian metric is
stratified, i.e. it is the disjoint union of smooth submanifolds of $V$. This
situation will be illustrated on the  case of complex Grassmann manifold.
Using a geometrical method, Wong\cite{wong,won,wo} has studied conjugate loci
 and cut loci of the Grassmann manifolds emphasising also their stratification.
Sakai\cite{sak1} has found out the cut locus of the  connected compact
symmetric manifold $V=U(n)/O(n)$, which  has $\pi_1(V)\cong{\gz}$.
 By refining the results of Ch. VII, \S 5
 "Control over singular set" from Helgason's book,\cite{helg}
 Sakai\cite{sak,sa}
studied the  cut locus of a point in a compact symmetric space  which is not
 necessarily simply connected and showed that it is
determined by the cut locus of a maximal totally geodesic flat torus of
 $V$. Takeuchi\cite{tak} has also proved the stratified structure of ${\bf CL}$
and ${\bf C}$ for compact symmetric manifolds. For other references see
 Kobayashi\cite{kob}. However, the expression of the conjugate locus as subset
of the Grassmann manifold is not known explicitely. This problem is largely
discussed elsewhere.\cite{viitor} In \S 4  of this paper we shall only collect
the main results on  this problem.

Most considerations in this Section concern only manifolds with the property

$${\rm Exp}\vert _o=\lambda \circ \exp \vert _{\got m} .
\leqno A) $$

Here
 ${\got g}={\got k}\oplus {\got m}$ is the orthogonal
decomposition with respect to the $B$-form as explained below at $B)$,
${\rm Exp}_p:\widetilde {\bf M}_p\rightarrow \widetilde {\bf M}$ is the
geodesic exponential map (cf. Ref. \onlinecite{helg} p. 33) and
 $\exp :{\got g} \rightarrow G $.

 In fact,
$A)$ expresses  that {\it the geodesics in} $\widetilde {\bf M}$
{\it are images of one-parameter subgroups of } $\widetilde {\bf M} \approx
G/K$. The symmetric spaces have property $A)$ (cf. Thm. 3.3 p. 208 in
Ref. \onlinecite{helg}).

We shall also be concerned with manifolds $\widetilde {\bf M}$
verifying the following condition:

{\it B)} On the Lie algebra {\got g} of $ G$
there exists an $ Ad(G)$-invariant,
symmetric, non-degenerate bilinear form $ B$ such that the restriction
of $B$ to the Lie algebra {\got k}  of $ K$  is likewise
non-degenerate.

We point out that {\it if the homogeneous space} $\widetilde {\bf M}\approx
G/K$ {\it verifies} $B)$, {\it then it also verifies} $A)$ (cf.
Corollary 2.5, Thm. 3.5 and Corollary 3.6 Chapter X in Ref. \onlinecite{kn}).
 Indeed, if
${\got g}={\got k}\oplus {\bf m}$ is the orthogonal
decomposition relative  to the $B$-form on {\got g}, then {\got m} is
canonically identified with the tangent space at $o,~
 \widetilde {\bf M}_o$. $  B) $ implies a (possibly indefinite)
$G$-invariant metric on $\widetilde {\bf M}$. It follows that $G/K$
is reductive, i.e. $[{\got k},{\got k}]\subset {\got k}$ and $[{\got
k},{\got m}]\subset {\got m}$. If $B)$ is true, then $\widetilde {\bf M}$
is naturally reductive (see p. 202 in Ref. \onlinecite{kn}) and $A)$ is also
verified.
The symmetric spaces verify besides the conditions of reductive
spaces, the condition $[{\got m},{\got m}]\subset {\got k}$ and, of course,
$A)$ is  verified too (see Thm. 3.2 Ch. Xl in Ref. \onlinecite{kn}).

Thimm\cite{tim}
 furnishes as another examples of homogeneous spaces verifying $B)$,
 besides the symmetric spaces, the Lie groups with bi-invariant metric
and the normal homogeneous spaces (i.e. $B$ is positive definite).
Kowalski\cite{kow}  studied generalised symmetric spaces still
 verifying condition
$A)$. See also Montgomery.\cite{mont}

2. Now we remember that in Ref. \onlinecite{sbl}
 we did the following Remark, which is in
 fact E. Cartan's theorem (see e.g. Thm. 3.3 p. 208  in Ref. \onlinecite{helg})
on geodesics on symmetric spaces expressed in the coherent state
setting:

\begin{rem}
 {\it The vector} $\vert
tB,j>=\exp\pi_j^{*\prime}(tB)\vert j>\in {\bf M} , B\in {\got m},$
{\it describes trajectories in} {\bf M} {\it corresponding to the image in the
manifold of coherent states} $\widetilde {\bf M}\hookrightarrow {\bf PL}$
{\it of geodesics through the identity coset element on the symmetric space}
$X\approx G/K$. {\it The dependence} $Z(t)=Z(tB)$ {\it appearing when
one passes from} eq. (\ref{b}) {\it to} eq. (\ref{z}) {\it describes in
${\cal V}_0$ a geodesic.}
\end{rem}
We shall reformulate  Remark 1 in a way very useful even for
practical calculations. The proof presented below, true in the particular
case of hermitian symmetric spaces, implies also Thm. 1.

\begin{rem}
 {\it For an} $n-$ {\it dimensional manifold} $X\approx G/K$
{\it which has Hermitian symmetric space structure, the parameters}
 $B_{\varphi}$
 {\it in formula} (\ref{b}) {\it of normalised coherent states are normal
coordinates
in the normal neighbourhood} ${\cal V}_0\approx {\gc}^n${\it around the point}
$Z_{\varphi}=0$ {\it on the manifold} $X$.
\end{rem}
{\it Proof}. The Harish-Chandra embedding theorem can be used (cf. e.g.
Ref. \onlinecite{knapp}; see also Ref. \onlinecite{sbl} for the present
 context). This theorem asserts that the map
$M^+\times K^{\gcm}\times M^- \rightarrow G^{\gcm}$ given by
 $(m^+,k,m^-)\rightarrow
m^+km^-$ is a complex analytic diffeomorphism onto an open dense subset
of $G^{\gcm}$ that contains $G_n$. Let ${\got m}^{\pm}$ be the $\pm i$
 eigenspaces
of $J$ and  $M^{\pm}$  the (unipotent, Abelian) subgroups of $G^{\gcm}$
corresponding to ${\got m}^{\pm}$. Then, in particular,
 $b:{\got m}^+ \rightarrow X_c=G^{\gcm}/P,~b(X)=\exp(X)P$
is a complex analytic diffeomorphism of ${\got
m}^+$ onto a dense subset of $X_c$ (that contains $X_n$) and the
Remark follows because the requirement $A)$ is fulfilled for the
symmetric spaces.

The proof of Remark 2 also shows that {\it the dependence} $Z=Z(B),$
with $ B\in {\got m}^+,$ and $ Z$
parametrizing  $\widetilde {\bf M}$, obtained passing
from eq. (\ref{b}) to (\ref{z}) using the relations (\ref{act})
 (the Baker-Campbell-Hausdorff formulas),\cite{sbl}
{\it expresses} in fact {\it the geodesic exponential} ${\rm Exp}_0:
 \widetilde {\bf M}_0
\rightarrow  \widetilde {\bf M}$. So we have proved the following \pagebreak
\begin{th}
 {\it Let}  $\widetilde {\bf M}$ {\it be a coherent state
manifold with Hermitian symmetric space structure,
parametrized in} ${\cal V}_0$ {\it around} $Z=0$ {\it as in eqs.} (\ref{z}),
(\ref{b}). {\it Then the
conjugate locus of the point} $ o $ {\it is obtained vanishing the Jacobian
of the exponential map $Z=Z(B)$ and the corresponding transformations of the
 chart from ${\cal V}_0$} .
\end{th}
The situation is very transparent in the case of the complex Grassmann
manifold $X_c=G_n({\gc}^{n+m})=SU(n+m)/S(U(n)\times U(m))$
and his noncompact dual $X_n=SU(n,m)/S(U(n)\times U(m))$. There\cite{sbl}

\begin{eqnarray}
\nonumber X_{n,c} & = & \exp\left(\matrix{0&B\cr
                            \pm B^*&0\cr}\right)o=
   \left(\matrix{{\rm co}\sqrt{BB^*}&B{\displaystyle
{ {\rm si}\sqrt{B^*B}\over \sqrt{B^*B}}}\cr
                                 \pm{\displaystyle {{\rm si} \sqrt{B^*B}\over
                         \sqrt{B^*B}}}B^*&{\rm co}\sqrt{B^*B}\cr}\right)o \\
\nonumber &  & \\
 & = & \left(\matrix{{\oo}&Z\cr 0&{\oo}\cr}\right)
\left(\matrix{({\oo}\mp ZZ^*)^{1/2}&0\cr 0&({\oo}\mp Z^*Z)^{1/2}\cr}\right)
\left(\matrix{{\oo}&0\cr \pm Z^*&{\oo}\cr}\right)o\\
\nonumber &  & \\
\nonumber & = & \exp\left(\matrix{0&Z\cr 0&0\cr}\right)P ,
\end{eqnarray}
where $B^*$ denotes the hermitian conjugate of the matrix $B$.
co is an abbreviation for the circular cosine cos (resp. the hyperbolic
cosine coh) for $X_c$ (resp. $X_n$) and similarly for si. The - (+) sign
in the equation above corresponds to the compact (resp. noncompact) $X$.

Here $Z$ and $B$ are $n\times m$ matrices related by the relation

\begin{equation}
Z=B{{\rm ta} \sqrt{B^*B}\over \sqrt{B^*B}} ,\label{geo}
\end{equation}
and ta is an abbreviation for the hyperbolic tangent tgh (resp. the
circular tangent tg) for $X_n$ (resp. $X_c$).
The dependence $Z=Z(B)$ describes in fact ${\rm Exp}:
G_n({\gc}^{n+m})_e \rightarrow G_n({\gc}^{n+m})$ in ${\cal V}_0$.
 Indeed, the equation of geodesics for $X_{c,n}$ is\cite{viitor}

\begin{equation}
\frac{d^2Z}{dt^2}-2\epsilon \frac{dZ}{dt}Z^+
(\oo +\epsilon ZZ^+)^{-1}\frac{dZ}{dt}=0~,\label{65}
\end{equation}
where $\epsilon =1~(-1)$ for $X_c$ (resp. $X_n$). It is easy to see that
 (\ref{geo})
verifies (\ref{65}) with the initial condition $\dot{Z}(0)=B$.

 $Z$ and $B$ in the  eq. (\ref{geo})  of geodesics are in the same time
 the parameters describing
 the coherent states in the paramerization given by eq. (\ref{z}) and
 respectively (\ref{b}).

3. Firstly, let us introduce a notation for the {\it polar divisor}
 of $\vert 0>\in {\bf M}:$
\begin{equation}
\Sigma _0=\left\{ \vert \psi > \vert \vert\psi >\in {\bf M},<0\vert \psi>=0
\right\} .
\end{equation}

This denomination is inspired from that one used by Wu\cite{wu}
 in the case of
the Grassmann manifold.

We shall prove the following \newpage
\begin{th}
 {\it Let} $\widetilde {\bf M}$ {\it be a
 homogeneous manifold}
$\widetilde {\bf M}\approx G/K$. {\it Suppose that there exists an unitary
irreducible representation} $\pi _j$ {\it of} $G$ {\it such that in a
neighbourhood }$ {\cal V}_0$ {\it around} $Z=0$ {\it the coherent states are
parametrized as in eq.} (\ref{z}). {\it Then the manifold} $\widetilde {\bf M}$
{\it can
be represented as the disjoint union}

\begin{equation}
\widetilde {\bf M} ={\cal V}_0\cup \Sigma _0.\label{reu}
\end{equation}

{\it Moreover, if the condition} $B)$ {\it is true, then}

\begin{equation}
\Sigma _0={\bf CL}_0,
\end{equation}
{\it and}
\begin{equation}
{\cal V}_0={\bf I}_0.
\end{equation}
\end{th}
{\it Proof.} We can take $\vert \psi >~=\vert \psi (Z)>~\in {\bf M}$
 such that
the parameters $Z$ are in ${\gc}^n$ as in formula (\ref{z}). Now, the
second relation (\ref{act}) implies that $<0 \vert \psi >~=1$ for
 $\vert \psi >\in
\xi^{-1}_0 ({\cal V}_0)$. It follows that the equation

\begin{equation}
 \cos \theta =0,
\end{equation}
where

\begin{equation}
\cos\theta ={\vert <0\vert \psi >\vert\over \Vert 0\Vert ^{1/2}\Vert
\psi \Vert ^{1/2}}=\Vert\psi\Vert ^{-1/2}, \label{unghi}
\end{equation}
does not have solutions for $\vert\psi >\in \xi^{-1}_0({\cal V}_0) $,
and the representation (3.4) follows.

To prove relation (3.5) if $B)$ is true, use is made of Thm. 7.4 and the
subsequent remark at p. 100 from Ref. \onlinecite{kn}, reproduced at the
 beginning of this section
in an enriched version. The theorem essentially says that any
 Riemannian manifold $\widetilde {\bf M}$ is the disjoint union of the
cut locus (closed cell) and the largest open cell of $\widetilde {\bf M}$
on which normal coordinates can be defined. But $Z\in {\gc}^n$ for
points of ${\cal V}_0$ corresponding to the largest normal coordinates
$B\in {\got m}$, because $B)$ implies $A)$.

Farther we shall prove a Corollary of Thm. 2. This is related to the
angle $\theta$ appearing in eq. (\ref{unghi}).

Firstly, let us introduce the (hermitian elliptic)  Cayley distance\cite{cay}
 in the projective space.
Let $(\cdot ,\cdot )$ be the scalar product in ${\bf K}$. If $\xi :
{\bf K}\setminus \{ 0\}\rightarrow {\bf PK}$ is the natural projection
$~ \xi : \omega \rightarrow [\omega ]$, then the Cayley distance is

\begin{equation}
d_c([\omega '],[\omega ])=\arccos \frac{\vert(\omega ',\omega)\vert}
{\Vert \omega '\Vert \Vert \omega \Vert }~.
\end{equation}
The infinite dimensional case is argued in Ref. \onlinecite{kobi}.
Before proving the Corollary, we shall present\newpage

\begin{rem} {\bf (Geometrical significance of transition amplitudes for
coherent states)} {\it Let }
 $\vert {\underline Z} >\in {\bf M}, Z\in {\cal V}_0$
{\it as in} (\ref{z}) {\it and} $\iota :\widetilde
 {\bf M}\hookrightarrow {\bf PL}$
 {\it the
embedding of the coherent state manifold into the  projective space.
Then the angle} $\theta =\theta (Z,Z')$ {\it defined by}
\begin{equation}
\theta \equiv \arccos \vert <\underline {Z'}\vert \underline {Z}>\vert
\end{equation}
{\it
is equal with the geodesic distance joining} $\iota (Z)${\it and} $\iota (Z')$,
\begin{equation}
\theta = d_c(\iota (Z'),\iota (Z)).
\end{equation}
{\it More generally, the (Cauchy) formula is true}
\begin{equation}
<\underline {Z'}\vert \underline {Z}>
=\frac{(\iota (Z'),\iota (Z))}{||\iota (Z')||\,||\iota (Z)||}.\label{cauchy}
\end{equation}
\end{rem}
{\it Proof}. The relation (\ref{cauchy}) is an immediate consequence of the
 fact
 that the complex analytic line bundle {\bf M} over $\widetilde {\bf M} $
is {\it  projectively induced} (see p. 139 in Ref. \onlinecite{hirz}), i.e.
 the coherent state manifold {\bf M} is the pull-back of the
 hyperplane bundle $H=[1]$ on ${\bf PL}$, i.e.
${\bf M}=\iota^*[1]$.\cite{cgr}

The denomination of eq. (\ref{cauchy}) as the Cauchy formula is due to the
 fact that
for the Pl\"ucker embedding of the Grassmann manifold this formula
is nothing else than the (Binet-) Cauchy formula.\cite{gant}

Combining Remark 3 and Theorem 2, we get the following
\begin{cor}
 {\it Suppose that} $\widetilde {\bf M}$
{\it is an homogeneous
manifold verifying} $ B)$ {\it and admitting the embedding}
$\iota :\widetilde {\bf M}\hookrightarrow {\bf PL}$. {\it Let} $0, Z\in
\widetilde {\bf M}$. {\it Then} $Z\in {\bf CL}_0$ {\it iff the Cayley distance
between the images } $\iota (0), \iota (Z)\in {\bf PL}$ {\it is} $\pi /2$
\begin{equation}
d_c(\iota (0),\iota (Z))=\pi /2.
\end{equation}
\end{cor}

\section{\hspace{-4mm}.\hspace{2mm}  AN EXAMPLE: THE COMPLEX GRASSMANN
 MANIFOLD }

\hspace*{0.8cm} The results of Section 3 will be illustrated on the
example of the complex Grassmann manifold. The calculation of the cut locus
on $G_n({\gc}^{n+m})$ was announced by Wong\cite{wong}
 and now  more
proofs (see e.g. Sakai\cite{sak} and also Ref. \onlinecite{kobi}) are
 available.
 Also Wong\cite{won}
 has announced the conjugate locus on
the Grassmann manifold, but, as far as I know, the proof has not been
published. Even more, the results of Wong
on conjugate locus on Grassmann manifold were contested by Sakai,\cite{sak}
who showed that the result of Wong is incomplete.

 The explicit calculation of the conjugate locus in the manifold
using Theorem 1 is presented elsewhere.\cite{viitor} Another proof of
 the results of Sakai referring to the tangent
conjugate locus is also presented there. Here we just indicate the
 parameters
appearing in the calculation in order to illustrate how the assertions of
 Section 3
referring to the cut locus and conjugate locus work in an concrete example.
 However, we do not have an explicit expression
for the part of the conjugate locus lost by Wong and only a geometrical
 characterisation in terms of the stationary angles.

1. Firstly we fix the notation concerning the geometric construction of
coherent state manifold when
$\widetilde {\bf M}$ is the complex Grassmann manifold (the manifold of
Slater determinants\cite{csb}).

Let {\bf O} be the $n-$plane passing through the origin of
 ${\gc}^N (N=n+m)$
corresponding to $Z=0$ in ${\cal V}_0\subset G_n({\gc}^N)$ in the
 representation (\ref{z}). Then $Z\in {\cal V}_0\approx {\gc}^{n\times m}$
 iff there
are $n$ vectors $\mb{z}_1,\ldots ,\mb{z}_n \in {\gc}^N$ such that
\begin{equation}
Z=\mb{z}_1\wedge\ldots\wedge\mb{z}_n\not=0 ~.\label{zzz}
\end{equation}
 We use the Pontrjagin
coordinates. Fixing the canonical basis
$\mb{e}_1,\ldots ,\mb{e}_N$ for $ {\gc}^N$, then

\begin{equation}
\mb{z}_i=\mb{e}_i+ \sum_{\alpha =n+1}^N Z_{i\alpha}\mb{e}_{\alpha},
\end{equation}

If the weight $j$ is taken as\cite{sbl}
\begin{equation}
j=(\underbrace{1,\ldots ,1}_n,\underbrace{0,\ldots ,0}_m),
\end{equation}
then we have the {\it equality\cite{viitor} of the scalar product}
 $<\cdot \vert \cdot >$ {\it of coherent vectors from } {\bf M} {\it and
of the hermitian scalar product} $((\cdot ,\cdot ))$
{\it in the holomorphic line bundle} $ \det ^*$\cite{chern}
\begin{equation}
<Z'\vert Z>=((\hat{Z}',\hat{Z}))=
\det ((\mb{z'}_i,\mb{z}_j))_{i,j=1,\ldots ,n}=
\det ({\oo}_n+ZZ'^*) ~.
\end{equation}
We have  used the notation
\begin{equation}
\hat{Z}= ({\oo}_n,Z),
\end{equation}
where $Z$ is an $n\times m$ matrix and ${\oo}_n$ is the unity
 $n\times n$ matrix.

So, {\it the parameters} $Z$ {\it in formula} (\ref{z}) {\it for
the Grassmann manifold of
coherent states are the Pontrjagin}\cite{pont} {\it coordinates}
$Z$ {\it in formula} (\ref{zz}).

Let us also introduce the Pl\"ucker coordinates $Z^{i_1\ldots i_n}$, i.e.
\begin{equation}
Z=\sum_{1\leq i_1<\ldots <i_n\leq N} Z^{i_1\ldots i_n}
\mb{e}_{i_1}\wedge\ldots\wedge\mb{e}_{i_n}	~.
\end{equation}

 Let $\iota: G_n({\bf K})
\hookrightarrow {\bf PL}$ be the Pl\"ucker embedding, where ${\bf K}=
{\gc}^N$, ${\bf L}={\gc}^{*N(m)}$, $N(m)={N\choose n}-1$. Using the
notation of Section 2 for $X=X_c=G_n({\gc}^N)$, then ${\bf M}'={\bf M}$,
that is $m_0=1$ in ${\bf M}'^{m_0}={\bf M}$, (i.e. the line bundle $\det ^*$
is not only ample, but very ample\cite{ss}) and ${\bf M}=\iota^*[1]$,
 where [1] is
the hyperplane section $H$ in {\bf L}.

The (Binet-) Cauchy formula\cite{gant} invoked in eq. (\ref{cauchy}) reads
 explicitly
\begin{equation}
\det ((\mb{z'}_i,\mb{z}_j))_{i,j=1,\ldots ,n}=
\sum_{1\leq i_1<\ldots <i_n\leq N} Z^{i_1\ldots i_n}
\overline {Z}'^{i_1\ldots i_n}~.
\end{equation}

2. Now we fix the notation referring to the Schubert varieties.

Let the sequences of integers
\begin{equation}
\omega=\{ 0\leq\omega (1)\leq \ldots \leq\omega (n)\leq m\},
\end{equation}
\begin{equation}
\sigma (i)=\omega (i)+i,~ i=1,\ldots ,n.
\end{equation}

The Schubert varieties  are defined as\cite{pont}
\begin{equation}
Z(\omega )=\left\{ X\in G_n({\gc}^{n+m}) \vert
 \dim (X\cap {\gc}^{\sigma (i)})\geq i\right\} .
\end{equation}

$Z(\omega ) $ are closed cells in the Grassmann manifold. The "jumps"
sequence\cite{mil} is introduced as
\begin{equation}
I_{\omega }=\left \{ i_0<i_1< \ldots <i_{l-1}<i_{l}=n\right\}~,
\end{equation}
where
\begin{equation}
\omega (i_h)<\omega (i_{h+1}), \omega (i)=\omega (i_{h-1}), i_{h-1}<i
\leq i_h, h=1,\ldots ,l.
\end{equation}

Let us consider the subset of  generic elements of $Z(\omega )$\cite{pont}
\begin{equation}
Z'(\omega )=\left\{ X\subset  G_n({\gc}^{n+m}) \vert
\dim (X\cap {\gc}^{\sigma (i_h)})= i_h,~  i_h\in I_{\omega}\right\} .
\end{equation}

The condition to get generic elements $Z$ of $Z(\omega )$,
$Z\in {\it V}_0\cap
Z(\omega )\subset Z'(\omega )$, is :\cite{pont,eh}

\begin{equation}
Z_{ij}=0,~~ j>\omega (i),~ i=1,\ldots ,n .
\end{equation}

Let also the notation
\begin{equation}
V^p_l=\left \{Z\subset  G_n({\gc}^{n+m})\vert \dim (Z\cap {\gc}^p)
\geq l\right\} ,
\end{equation}
\begin{equation}
W^p_l=V^p_l-V^p_{l+1}=
\left \{Z\subset  G_n({\gc}^{n+m})\vert \dim (Z\cap {\gc}^p)
= l\right\} ,
\end{equation}
\begin{equation}
\omega^p_l=(\underbrace{p-l,\ldots ,p-l}_l, \underbrace{m,\ldots ,m}_{n-l}).
\end{equation}

Then\cite{won,wo,viitor}
\begin{equation}
V^p_l=Z(\omega ^p_l);~~ W^p_l=Z'(\omega ^p_l)~ .
\end{equation}

3. Now we briefly remember some notions referring to the stationary angles.

  Let $Z', Z$ be two $n-$planes of $G_n({\gc}^{n+m})$
given as in eq. (\ref{zzz}). Then the ($n$) {\it stationary angles}
(see Jordan\cite{jor} for the real case),
{\it of which most} $r=\min (m,n)$ {\it are nonzero}, are defined as the
{\it stationary} angles $\theta\in[0,\pi /2]$ between the vectors
\begin{equation}
\mb{a}=\sum_{i=1}^na_i\mb{z'_i},~\mb{b}=\sum_{i=1}^nb_i\mb{z}_i,~
\end{equation}
where
\begin{equation}
\cos \theta =\frac{|(\mb{a},\mb{b})|}{||\mb{a}||||\mb{b}||}~.\label{unghii}
\end{equation}
We remember the following two Lemmas\cite{jor,ros,rosi,viitor}

\begin{lem}\label{lem4}
 The squares
 $\cos ^2\theta _i$ {\it of the stationary angles between the\/ }
$\mbox{}n-${\it planes} $Z, Z'$ {\it with} $((Z,Z'))
\not= 0 $
 are given
as the eigenvalues of a matrix $W$  which, for $Z, Z'\in {\cal V}_0$
is

\begin{equation}
W=({\oo}+ZZ^+)^{-1}({\oo}+ZZ'^+)({\oo}+Z'Z'^+)^{-1}
({\oo}+Z'Z^+)~ .\label{53}
\end{equation}
\end{lem}
\begin{lem}\label{lem5}
  Let $\theta$  be the angle defined by the
hermitian scalar product in  the following equation

\begin{equation}
\cos \theta (Z',Z)\equiv\frac{|((Z',Z))|}{\Vert Z'\Vert \Vert Z
\Vert }~= {|\det ({\oo}+ZZ'^+) |
 \over |\det ({\oo}+ZZ^+)|^{1/2} |\det ({\oo}+Z'Z'^+)|^{1/2}} ~,\label{10.1}
\end{equation}
 $d_c$  the Cayley distance
and  $\theta _1,\ldots ,\theta_n$  the stationary angles. Then
\begin{equation}
\cos \theta (Z,Z')= \cos d_c(\iota (Z'),\iota (Z))=
\cos \theta_1\cdots\cos \theta _n~.
\label{55}
\end{equation}
\end{lem}
 It can be
proved\cite{viitor}
 that the eigenvalues of W appear also in the expression of
the distance on the complex
Grassmann manifold (see also Siegel\cite{sieg}).

Note also that if the expression (\ref{geo}) of the dependence $Z=Z(B)$ is
introduced in the formula of the distance between the points
 $Z=0$ and $Z\in {\it V}_0$ on the Grassmann
manifold, then
\begin{equation}
d^2=\sum |B_{ij}|^2 .
\end{equation}

The last equation expresses the fact that the parameters $B$ in
eq. (\ref{b})
of coherent states are indeed the normal coordinates as it is
asserted in Remark 2.

4. Below we present the cut locus and the conjugate locus for \Gras .

 ${\bf O}^{\perp}$ denotes the orthogonal complement of the $n$-plane
${\bf O}$ in ${\gc}^N$.

\begin{rem}[Wong\cite{wong}] {\it The cut locus of the point }{\bf O}
 {\it is given by}
\begin{eqnarray}
\nonumber {\bf CL}_0 & = & \Sigma_0  =  V^m_1 = Z(\omega^m_1)= Z(m-1,m,\ldots ,
m)\\ & =  & \left\{ X\subset G_n({\gc}^{n+m})
\vert \dim (X\cap {\bf O}^{\perp})\geq 1\right\} .
\end{eqnarray}

{\it The cut locus in \Gras~ is given by those $n-$planes which have
at least one of the stationary angles $\pi /2$ with the plane {\bf O}}.
\end{rem}

{\it Proof.} An immediate proof can be obtained using the results of
Wu
referring to the polar divisor $\Sigma_0$ on the Grassmann
manifold  (see Ch. l in  Ref. \onlinecite{wu}) and the
theorems  characterising the canonical (universal, det)
bundle on $G_n({\gc}^N)$ (see especially Prop. 3.3 Ch. 7 in
Ref. \onlinecite{hus}), which are particularisations of the representation
in Thm. 2.

The following theorem summarize the known facts about the tangent conjugate
locus and conjugate locus in \Gras .\cite{wong,sak,viitor} The relevant fact
for the  present paper is that the conjugate locus can be calculated using
Theorem 1.
\begin{th}
\label{sakth}
The tangent conjugate locus $C_0$ of the point ${\bf O}\in\Gras$ is given by
\begin{equation}
\label{ura}
C_0=\bigcup_{k,p,q,i}ad\,k(t_iH)~,~i=1,2,3;~1\leq p<q\leq r,
\end{equation}
where the vector $H\in{\got a}$
   is normalised,
\begin{equation}
H=\sum_{i=1}^r h_iD_{i\,n+i},~h_i\in\gr,~\sum h^2_i=1~.
\label{hh}
\end{equation}
 The parameters $t_i,~i=1,2,3$ in eq. (\ref{ura}) are
\begin{equation}
\begin{array}{l@{\:=\:}c@{\:,\:}l}
t_1  & \displaystyle{\frac{\lambda \pi}{|h_p\pm h_q|}}  &
 ~\mbox{\rm multiplicity}~2;\\[2.ex]
t_2  & \displaystyle{\frac{\lambda \pi}{2|h_p|}} &
  ~\mbox{\rm multiplicity}~1;\\[2.ex]
t_3  & \displaystyle{\frac{\lambda \pi}{|h_p|}}  &
 ~\mbox{\rm multiplicity}~2|m-n|; ~\lambda\in \gz^{\star}~.
\end{array}
\label{ttt}
\end{equation}

 The conjugate locus of {\bf O} in \Gras  ~is given by the union
\begin{equation}
\label{reun}
{\bf C}_0={\bf C}^W_0\cup {\bf C}^I_0.
\end{equation}

 The following relations are true

\begin{equation}
\label{ect1}
{\bf C}^I_0= \exp \bigcup_{k,p,q}Ad\,k(t_1H)~,
\end{equation}
\begin{equation}
\label{ect2}
{\bf C}^W_0= \exp \bigcup_{k,p}Ad\,k(t_2H)~,
\end{equation}
i.e. exponentiating the vectors of the type $t_1H$ we get the points
of ${\bf C}^I_0$ for which at least two of the stationary angles with {\bf O}
are equal, while the vectors of the type $t_2H$ are sent to the points of
 ${\bf C}^W_0$ for which at least one of the stationary angles with {\bf O}
is $0$ or $\pi /2$.

The ${\bf C}^W_0$ part of the conjugate locus is given by the disjoint union
\begin{equation}
{\bf C}^W_0=\label{72}
\cases {V^m_1\cup V^n_1,&$ n\leq m,$\cr
            V^m_1\cup V^n_{n-m+1},& $n>m,$ \cr}
\end{equation}
 where
\begin{equation}
V^m_1=\label{74.1}
\cases {{\gcp}^{m-1},& for $n=1 ,$\cr
              W^m_1\cup W^m_2\cup \ldots W^m_{r-1}\cup W^m_r,&$1<n ,$\cr}
\end{equation}
\begin{equation}
W^m_r=\cases {G_r({\gc}^{\max (m,n)}),&$n\not= m,$\cr
             {\bf O}^{\perp},&$n=m ,$\cr}
\end{equation}
\begin{equation}
V^n_1=\label{72.2}
\cases {W^n_1\cup \ldots \cup W^n_{r-1}\cup {\bf O},&$1<n\leq m,$\cr
             {\bf O}, &$n=1 ,$ \cr}
\end{equation}
\begin{equation}
V^n_{n-m+1}=W^n_{n-m+1}\cup W^n_{n-m+2}\cup\ldots\cup W^n_{n-1}\cup
{\bf O}~,~ n>m ~.\label{73.3}
\end{equation}
\end{th}

{\it Sketch of the Proof}. The tangent conjugate locus $C_0$  for \Gras~  in
the case $n\leq m$ was obtained by
 Sakai.\cite{sak} Sakai has observed that Wong's result on the conjugate locus
in the manifold is incomplete, i.e. ${\bf C}^W_0\subset{\bf C}_0$ but
${\bf C}^W_0\ine {\bf C}_0=\exp C_0$. The proof of Sakai consists in
solving the eigenvalue equation $R(X,Y^i)X=e_iY^i$
which appears when solving  the Jacobi equation,
 where the
curvature for the symmetric space $X_c=G_c/K$ at $o$ is simply $R(X,Y)Z=
[[X,Y],Z],~X,Y,Z\in\got{m}_c$. Then $q={\rm Exp}_0tX$ is conjugate to $o$ if
$t=\pi \lambda/\sqrt{e_i},~\lambda\in\gz^{\star}\equiv \gz\setminus \{0\}$.

Above {\got a} is the Cartan subalgebra of the
symmetric pair $(SU(n+m),S(U(n)\times U(m)))$\cite{helg,sak,viitor}
consisting of vectors of the form (\ref{hh})
where $r$ is the symmetric rank of $X_c$ (and $X_n$) and we use the notation
$D_{ij}=E_{ij}-E_{ji},~i,j=1,\ldots ,N.$
$E_{ij}$ is the matrix with entry $1$ on line $i$ and column $j$ and $0$
otherwise.
 The results in the complex Grassmann manifold are obtained farther using the
exponential map given by eq. (\ref{geo}).

The same result on the calculation of the tangent conjugate locus can be
obtained\cite{viitor} using Prop. 3.1 p. 294 in the book of
 Helgason.\cite{helg}  This Proposition asserts that  $H\in
{\got a}$ is conjugate with $o$  iff $\alpha (H)\in i\pi\gz^{\star}$
 for some root $\alpha$  which do not vanishing identically on {\got a}. The
 eigenvalues
of the  equation $[H,X]=\lambda X,~\forall H\in {\got a}$,
lead\cite{viitor} to the values given  in equation (\ref{ura}) for the
 parameters $t_1-t_3$.

 The direct proof\cite{viitor}  in the Grassmann manifold
uses in Theorem 1 the dependence $Z=Z(B)$ furnished by eq.
 (\ref{geo})
which gives the geodesics on
$G_n({\gc}^{n+m})$  and the Jordan's
 stationary angles between two $n-$planes.
The stationary angles between two $n-$planes are given by Lemma 1
 and appear in the relation given by Lemma 2.

The proof\cite{viitor} is done in four steps. a) Firstly, a diagonalization of
 the $n\times
m$ matrix $Z$ is performed. b) Secondly, the Jacobian of a transformation
of complex dimension one is computed. c) The cut locus is reobtained and his
 contribution to the  conjugate locus is taken into account. d) The
nonzero angles are counted using the following property of the stationary
 angles: if the $n'~(n)$-plane (resp.  $Z_{n}$) are such that $Z'_{n'}
\cap Z_n=Z"_{n"}$, than $n'-n"$ angles of $Z'_{n'}$ and $Z_n$ are different
from $0$ and $n"$ are 0.

\section{\hspace{-4mm}.\hspace{2mm}  CONCLUSION AND DISCUSSION }
\hspace*{0.8cm}

In this paper it was shown that for a certain class of homogeneous manifolds
which include the symmetric ones there is a relationship between geodesics and
 coherent states. The starting point\cite{sbl} of  the present investigation,
contained in Remark 1, is
the observation that for symmetric spaces, if one expresses the
 parameters
 $Z$ in eq. (\ref{z}) as a function
of the parameters $B$  in eq. (\ref{b}), both characterising
 the coherent states, explicit local formulas  for the geodesic exponential map
are obtained. For Hermitian symmetric spaces the dependence $Z=Z(B)$
can be found using the Harish-Chandra decomposition or the so called
Baker-Campbell-Hausdorff formulas.\cite{sbl}  Thus Theorem 1 permits a
 calculation of the conjugate locus in
\Gras .   However, the explicit form of the conjugate locus in \Gras ~
is not completely known.\cite{viitor} The part of the conjugate locus
 ${\bf C}^W_0$
 determined by Wong is expressible as Schubert varieties,\cite{won} while the
 rest\cite{sak} ${\bf C}^I_0$ can be characterised\cite{viitor} as the subset
 of points of \Gras~ which have at least two of the stationary
 angles with the fixed $n$-plane ${\bf O}$ equal.   ${\bf C}^I_0$ contains as
 subset the maximal set of mutually isoclinic\cite{wwong} subspaces of
 the Grassmann
manifold, which are isoclinic spheres,\cite{wwong,ww1} with dimension given
 by the solution of the Hurwitz\cite{hur} problem. This part referring to the
 explicit calculation of the conjugate locus on \Gras \ was only briefly
 included in
 Sec. IV, the full details being presented elsewhere.\cite{viitor}

 The main remark of this paper contained in Theorem 2, the equality
 (\ref{basical}), is a simple consequence of the fact that
any manifold is the disjoint union of a maximal normal  neighbourhood
${\cal V}_0$ of a
 point $0$  and the cut locus ${\bf CL}_0$. It would be interesting to find a
geometrical description of the polar divisor for manifolds which are not
characterised by  {\it condition B)}.  On the other side, the problem to find
 explicitly
the cut locus    on nonsymmetric spaces is a difficult one.\cite{kob}
Also it was proved that for homogeneous manifolds verifying {\it
the condition B)} and admitting an embedding in an adequate projective Hilbert
 space a necessary and sufficient condition that a point to belong to the cut
 locus of another point
is that the Cayley distance between the images  of the
 points  through the embedding to be $\pi /2$. \nopagebreak This category of
 manifolds includes all the coherent states manifolds which admit
 prequantization.\cite{sb}

\vspace*{0.8cm}
\hspace{-4mm}\hspace{2mm}{\bf  ACKNOWLEDGEMENTS}
\hspace*{0.8cm}

\vspace{1cm}

 The author expresses his gratitude
to Prof. Anne Boutet de Monvel and CNRS for the opportunity to work during the
summers  1993 and 1994 in the
\'Equipe de Physique Math\'ematique et G\'eom\'etrie  at the Institute
 de Math\'ema\-tique de l'Universit\'e Paris 7 Denis Diderot.
The constant interest of Professor L. Boutet de Monvel
is kindly acknowledged. Discussions during the XIII and XIV Workshops on
Geometrical
 Methods in Physics at  Bia\l owie\.za, especially with
Professors D. Simms, M. Cahen, A. M. Perelomov and J. Klauder are
 acknowledged. The
author also expresses his thanks to Professors K. Teleman, S. Kobayashi,
 M. Berger, T. Sakai and Th. Hangan for  suggestions.



\begin{references}


\bibitem{klauder}{\it Coherent
States}, edited by J. R. Klauder and  B. S. Skagerstam
 (Word Scientific, Singapore, 1985).

\bibitem{sbaa}S. Berceanu, ``From quantum mechanics to classical mechanics
 and back, via coherent states'',  in
{\it Quantization and Infinite-Dimensional Systems} (Plenum, New York,
  1994),  p. 155.

\bibitem{per}A. M. Perelomov, ``Coherent states for arbitrary Lie groups'',
 Commun. Math. Phys. {\bf 26}, 222 (1972).

\bibitem{raw}J. R. Rawnsley, ``Coherent states and K\"ahler manifolds'',
Quart. J. Math. Oxford  {\bf 28}, 403 (1977).

\bibitem{kost}B. Kostant, ``Quantization and unitary representations'',
  in Lecture Notes in Mathematics, Vol. 170,
edited by C. T. Taam (Springer-Verlag, Berlin 1970), p. 87.

\bibitem{chern}S. S. Chern, {\it Complex Manifolds without Potential
Theory} (Van Nostrand, Princeton, 1967).

\bibitem{hus}D. Husemoller, {\it Fibre Bundles} (Mc Graw-Hill, New York
1966).

\bibitem{hirz}F. Hirzebruch, {\it Topological Methods in Algebraic Geometry}
(Springer-Verlag, Berlin, 1966).

\bibitem{cval}V. Ceausescu and A. Gheorghe, ``Classical limit and
quantization of Hamiltonian systems'',  in  {\it Symmetries and Semiclassical
Features of Nuclear Dynamics}, Lecture Notes in Physics, Vol.
279, edited by A. A. Raduta (Springer-Verlag,  Berlin, 1987), p. 69.

 \bibitem{sbl}S. Berceanu and L. Boutet de Monvel, ``Linear dynamical
systems, coherent state manifolds, flows and matrix Riccati equation'',
 J. Math. Phys.
 {\bf 34}, 2353 (1993).

\bibitem{kn}S. Kobayashi and K. Nomizu, {\it Foundations of Differential
Geometry}, Vol. ll (Interscience, New York, 1969).

\bibitem{helg}S. Helgason, {\it Differential Geometry, Lie groups and
Symmetric Spaces} (Academic, New York, 1978).

\bibitem{wong}Y. -C. Wong, ``Differential Geometry of Grassmann manifolds'',
 Proc. Nat. Acad. Sci. U.S.A. {\bf 57}, 589 (1967).

\bibitem{won}Y. -C. Wong, ``Conjugate loci in Grassmann manifold'',
Bull. Am. Math. Soc. {\bf 74}, 240 (1968).

\bibitem{sak}T. Sakai, `` On cut loci on compact symmetric spaces'',
 Hokkaido Math. J. {\bf 6}, 136 (1977).

\bibitem{kob}S. Kobayashi, ``On conjugate and cut loci'',  in
{\it Global Differential Geometry}, M.A.A. Studies in Mathematics,
Vol.  27, S. S. Chern editor  (1989), p. 140.

\bibitem{viitor}S. Berceanu, ``On the Geometry of complex Grassmann manifold,
its noncompact dual and coherent states'', preprint Bucharest, Institute
of Atomic Physics FT-409-1995, September.

\bibitem{jor}C. Jordan, ``Essai sur la G\'eom\' etrie \`a n dimensions'',
 Bull. Soc. Math. France {\bf t. lll}, 103 (1875).

\bibitem{sb}S. Berceanu, ``The coherent states: old geometrical methods
in new quantum clothes'', poster presented at the XI$^{\rm {\grave{e}me}}$
 Congr\`es
International de Physique Math\'ematique, Paris, 18-23 Juillet (1994);
 preprint Bucharest, Institute of Atomic Physics, FT-398-1994 and
 preprint Universit\"{a}t  Bielefeld BiBoS Nr. 664/11/94.


\bibitem{sbpol}S. Berceanu, ``Coherent states and global Differential
  Geometry'', talk presented at the XIII Workshop on Geometrical Methods in
 Physics at Bia\l owie\.za,  Poland, July 9-15 (1994), to appear in the
Proceeding of the Workshop {\it Quantisation, coherent states, and complex
 structures}, edited by J. P. Antoine, S. T. Ali, W. Lisiecki, I. Mladenov and
A. Odzijewicz, Plenum (1995).

\bibitem{gs}V. Guillemin and S. Sternberg, {\it Symplectic Techniques in
Physics} (Cambridge U. P., London, 1984).


\bibitem{od}A. Odzijewicz, ``Coherent states and geometric quantization'',
  Commun. Math. Phys.  {\bf 150}, 85 (1992).

\bibitem{kir}A. A. Kirillov, {\it Elements of Theory of Representations}
(Springer-Verlag, New York, 1976).

\bibitem{on}E. Onofri,  ``On quantization theory for homogeneous
K\"ahler spaces'', preprint Parma IFPR-T-038 (1974).

\bibitem{wood}H. Woodhouse, {\it Geometric Quantization} (Oxford U. P.,
 Oxford, 1980).

\bibitem{knapp}A. W. Knapp, {\it Representation Theory of Semisimple Lie
Groups} (Princeton, NJ, 1986).

\bibitem{cgr}M. Cahen, S. Gutt and J. Rawnsley, ``Quantization on K\"ahler
manifolds'', ll,  Trans. Math. Soc. {\bf 337}, 73 (1993).

\bibitem{ss}B. Shiffman and A. J. Sommese, {\it Vanishing Theorems on Complex
Manifolds}, Progress in Mathematics, Vol. 56 (Birkh\"auser, Boston,
1985).

\bibitem{gh}P. Griffith and J. Harris, {\it Principles of Algebraic Geometry}
(Wiley, New York, 1978).

\bibitem{sbcag}S. Berceanu and A. Gheorghe, ``On the construction of
perfect Morse functions on compact manifolds of coherent states'',
 J. Math. Phys.  {\bf 28}, 2899 (1987).

\bibitem{kobi}S. Kobayashi, ``On the Geometry of bounded domains'',
 Trans. Amer. Math. Soc. {\bf 92}, 267 (1959).

\bibitem{lis}W. Lisiecki, ``A classification of  coherent state
representations of unimodular Lie groups'',  Bull. Amer. Math.
Soc. {\bf 25}, 37 (1991); ``Coherent state representations. A survey'',
 to appear
in  Math. Rep. (1995).

\bibitem{cr}R. Crittenden, ``Minimum and conjugate points in symmetric
 spaces'',  Canad. J. Math. {\bf 14}, 320 (1962).

\bibitem{war}F. W. Warner, ``The conjugate locus of a Riemannian manifold'',
 Amer. J. Math. {\bf 87}, 575 (1965).

\bibitem{wei} A. D. Weinstein, ``The cut locus and conjugate locus of a
Riemannian manifold'',  Ann. of Math. {\bf 87}, 29 (1968).

\bibitem{wo}Y. -C. Wong, ``A class of Schubert varieties'', J. Diff.
Geom. {\bf 4}, 37 (1970).

\bibitem{sak1}T. Sakai, ``The manifold of the Lagrangean subspaces of a
symplectic vector space'',  J. Diff. Geom. {\bf 12}, 555 (1977).

\bibitem{sa}T. Sakai, ``On the structure of cut loci in compact riemannian
 symmetric spaces'',  Math. Ann. {\bf 235} 129 (1978).

\bibitem{tak}M. Takeuchi, ``On conjugate loci and cut loci of compact symmetric
 spaces '' I, II,  Tsukuba Math. J.   {\bf 2}, 35 (1977); {\bf 3}, 1 (1979).



\bibitem{tim}A. Thimm, ``Integrable geodesic flows on homogeneous spaces'',
 Ergod. Theory Dyn. Syst. {\bf 1}, 495 (1981).

\bibitem{kow}O. Kowalski, {\it Generalised Symmetric Spaces}, Lecture Notes
in Mathematics Vol. 805 (Springer-Verlag, Berlin, 1980).

\bibitem{mont}R. Montgomery, ``Isoholonomic problems and some
applications'',  Comm. Math. Phys. {\bf 128}, 565 (1990).

\bibitem{wu}H. H. Wu, {\it The Equidistribution Theory of Holomorphic Curves},
Annals of Maths. Studies Vol. 164 (Princeton, New Jersey, 1970).

\bibitem{cay}A. Cayley, ``A sixth memoir upon quantics'',
 Phil. Trans. Royal. Soc. London {\bf 149}, 61 (1859).

\bibitem{gant}F. P. Gantmacher, {\it Teoria Matritz}
(Nauka, Moskwa, 1966).



\bibitem{csb}S. Berceanu and A. Gheorghe, ``Perfect Morse functions on
manifold of Slater determinants'',  Rev. Roum. Phys.  {\bf 34}, 125 (1989).

\bibitem{pont}L. C. Pontrjagin, ``Charakteristiceskie tzikly differentziruemyh
mnogobrazia'',   Mat. sb. {\bf 21}, 233 (1947).

\bibitem{mil}J. Milnor and J. D. Stasheff, {\it Characteristic Classes},
Annals of Maths. Studies  Vol. 76 (Princeton, New Jersey, 1974).

\bibitem{eh}C. Ehresman, ``Sur la topologie de certain espaces homogen\`es'',
  Ann. Math.  {\bf 35}, 396 (1934).




\bibitem{ros}B. Rosenfel'd, ``Vnnutrenyaya geometriya mnojestva
m-mernyh ploskastei n-mernova ellipticeskovo prostranstva'',
 Izv. Akad. Nauk. SSSR, ser. mat. {\bf 5}, 353 (1941).

\bibitem{rosi}B. Rosenfel'd, {\it Mnogomernye Prostranstva} (Nauka, Moskwa,
1966); {\it Neevklidovy Pronstranstva} (Nauka, Moskwa, 1969).

\bibitem{sieg}C. L. Siegel, {\it Symplectic Geometry} (Academic Press,
 New York, 1964).

\bibitem{wwong}Y. -C. Wong, ``Isoclinic $n$-planes in Euclidean $2n$-space,
Clifford parallels in elliptic $(2n-1)$ space, and the Hurwitz matrix
 equations'',
 Mem. Amer. Math. Soc. {\bf 41} (1961).

\bibitem{ww1}J. A. Wolf, ``Geodesic spheres in Grassmann manifolds'',  Ilin. J.
Math. {\bf 7}, 425 (1963); ``Elliptic spaces in Grassmann manifolds'',
  Ilin. J. Math. {\bf 7}, 447 (1963).


\bibitem{hur}A. Hurwitz, ``\"{U}ber die Komposition der quadratischen Formen
 von beliebig vielen Variablen'',  Nach. v. der Ges. der. Wiss., G\"{o}tingen
(Math. Phys. K1) {\bf 309}, (1898); reprinted in Math. Werke {\bf Bd. 2}, p.
 641.



\end{references}
\end{document}